\documentclass[letterpaper,12pt,leqno]{amsart}
\usepackage{fullpage}
\usepackage{amsfonts,amscd}
\usepackage[T1]{fontenc} 

\usepackage{mathptmx} 

\usepackage{tikz}
\usepackage{pgfplots}
\usepackage{microtype} 
\tikzset{>=latex} 
\usepackage{hyperref}

\usepackage[english]{babel}

\theoremstyle{plain}
\newtheorem{theorem}                {Theorem}      [section]
\newtheorem{proposition}  [theorem]  {Proposition}

\theoremstyle{definition}
\newtheorem{example}      [theorem]  {Example}
\newtheorem{remark}       [theorem]  {Remark}
\newtheorem{definition}   [theorem]  {Definition}

\def \i{\mathrm i}

\def \R{{\mathbb R}}
\def \s{{\mathbb S}}
\def \z{{\mathbb Z}}

\DeclareMathOperator{\cst}{constant}
\DeclareMathOperator{\grad}{grad}
\DeclareMathOperator{\trace}{Trace}
\DeclareMathOperator{\Div}{div}

\title[Biharmonic and Biconservative  Submanifolds]{An Invitation to Biharmonic and Biconservative  Submanifolds}

\author{Stefano Montaldo}
\address{Universit\`a degli Studi di Cagliari\\
Dipartimento di Matematica e Informatica\\
Via Ospedale 72\\
09124 Cagliari, Italy}
\email{stefano.montaldo@unica.it}

\begin{document}



\thanks{The author is member of the Italian National Group G.N.S.A.G.A. of INdAM.. The work was partially supported by the Project {ISI-HOMOS} funded by Fondazione di Sardegna and by PNRR e.INS Ecosystem of Innovation for Next Generation Sardinia (CUP F53C22000430001, codice MUR ECS00000038)}

\maketitle
\begin{center}
{\em \small  To the memory of Francesco Mercuri}
\end{center}

\begin{abstract}
  This note is based on a lecture delivered by the author at the Second Conference on Differential Geometry, held in Fez in October 2024. It offers an accessible introduction to biharmonic and biconservative submanifolds, exploring the motivations for their study and highlighting some key facts and open problems in the field.
\end{abstract}

\section*{Aim of the Study}

Minimal immersed submanifolds $\varphi:M^m\hookrightarrow (N^n,h)$ of a Riemannian manifold play a fundamental role in Differential Geometry. They are characterized by the vanishing of the {\em mean curvature vector field} ${\mathbf H}=\trace B/m$, where $B$ denotes the {\em second fundamental form} of the isometric immersion $\varphi$. 

A particularly elegant way to describe minimal  submanifolds is by embedding their definition within the framework of {\em harmonic mappings} in the sense of Eells-Sampson \cite{MR164306}. More precisely, let $C^{\infty}(M,N)$ denote the space of smooth maps $\varphi : (M^m,g)\to(N^n,h)$ between two Riemannian manifolds, then a map  $\varphi\in
C^{\infty}(M,N)$ is called {\it harmonic} if it is a critical point of the
{\it energy} functional
\begin{equation}\label{energia}
E:C^{\infty}(M,N)\to\R, \quad E(\varphi)=\frac{1}{2}\int_{M}\,|d\varphi|^2\; v_g.
\end{equation}
 In particular, $\varphi$ is harmonic if it is a solution of the Euler-Lagrange system of equations associated to \eqref{energia}, i.e.,
\begin{equation}\label{harmonicityequation}
  - d^* d \varphi =   {\trace} \, \nabla d \varphi =0,
\end{equation}
where $d$ is the exterior differential operator while $d^*$ represents the codifferential operator, see \cite[pag. 7-8]{MR703510} for a formal definition. 
The left member of \eqref{harmonicityequation} is a vector field along the map $\varphi$ or, equivalently, a section of the pull-back bundle $\varphi^{-1} TN$: it is called {\em tension field} and denoted $\tau (\varphi)$. In local coordinates $(U,x^i)$ on $M^m$ and $(V,y^{\alpha})$ on $N^n$ with $\varphi(U)\subset V$, system \eqref{harmonicityequation}, which consists of second order semilinear elliptic partial diﬀerential equations, acquires the form
$$
-\,\Delta \varphi^{\gamma}\,+\, ^N\Gamma_{\alpha \beta}^\gamma \,\frac{\partial \varphi^\alpha}{\partial x^i}\,\frac{\partial \varphi^\beta}{\partial x^j}
\,g^{ij}\, =\,0\,, \quad \gamma=1,\ldots,\,n,
$$
where $\Delta$ is the Beltrami-Laplace operator on $(M,g)$ and the $^N\Gamma$'s are Christoffel symbols of $(N,h)$.

When $\varphi:M^m\hookrightarrow (N^n,h)$ is an isometric immersion, the tension field of $\varphi:(M^m, \varphi^*h)\to (N^n,h)$ simplifies to
\[
\tau(\varphi)=\trace\nabla d\varphi=\trace B=m {\mathbf H}\,.
\]
This immediately implies that {\em harmonic isometric immersions} are precisely  {\em minimal immersions}.
Thus, the study of minimal immersions naturally leads to the investigation of harmonic isometric immersions.

More interesting, the elegant perspective of viewing minimal immersions as harmonic isometric immersions provides a natural pathway to generalizing the notion of minimal immersions by exploring suitable extensions of the notion of harmonic maps.

To pursue this idea, we will present in this note two natural generalizations of minimal   submanifolds.
\section{Biharmonic   submanifolds}
In \cite{MR0216519}, Eells and Sampson, just one year after their celebrated paper on harmonic maps \cite{MR164306}, proposed to study the critical points of the following {\em higher order energy functional}
\begin{equation}\label{ES-energia}
E_r:C^{\infty}(M,N)\to\R, \quad     E_r(\varphi)=\frac{1}{2}\int_{M}\,|(d^*+d)^r (\varphi)|^2\;v_g\,\, .
\end{equation}
When $r=1$ we recover the energy functional \eqref{energia}, while, when $r=2$, the functional \eqref{ES-energia} becomes 
\begin{align}
\label{bienergy}
E_2(\varphi)=&\frac{1}{2}\int_M|(d^*+d)(d^*+d)(\varphi)|^2 \ v_g=  \frac{1}{2}\int_M|d^*d\varphi)|^2 \; v_g\nonumber\\
=&\frac{1}{2}\int_M|\tau(\phi)|^2 \; v_g\,,
\end{align}
and it is called the \emph{bienergy functional}. Critical points of \eqref{bienergy} are called \emph{biharmonic maps} and can be expressed by the vanishing of the \emph{bitension field} \(\tau_2(\varphi)\) as follows (see \cite{MR886529}):
\begin{align}
\label{bitension}
\tau_2(\varphi):=-\bar\Delta\tau(\varphi)-\trace R^N(d\varphi(\cdot), \tau(\varphi))d\varphi(\cdot)=0\,,
\end{align}
where $\bar\Delta=d^*d$ is the rough Laplacian defined on sections of the pull-back bundle \(\varphi^{-1}TN\) and $R^N$
is the curvature operator on $(N^n,h)$. The biharmonic condition \eqref{bitension} constitutes a fourth-order semi-linear elliptic system of partial differential equations which is, in its generality, challenging to solve.  Nevertheless, over the past three decades, biharmonic maps have garnered significant attention and there is a vast literature on the subject. For a comprehensive overview of the current research status on biharmonic maps, we refer, for example, to the book \cite{MR4265170}.

Now, we say that an immersed submanifold  $\varphi:M^m\hookrightarrow (N^n,h)$ is a {\em biharmonic submanifold} if the immersion $\varphi$ is biharmonic as a smooth map  $\varphi:(M^m, \varphi^*h)\to (N^n,h)$. Thus a submanifold is biharmonic if the mean curvature vector field $\mathbf H$ satisfies:
\begin{align}
\label{bitension-H}
\bar\Delta\mathbf H+\trace R^N(d\varphi(\cdot), \mathbf H)d\varphi(\cdot)=0.
\end{align}

Clearly, a minimal submanifold, that is, one with ${\mathbf H}=0$, is biharmonic. Consequently, the class of biharmonic submanifolds extends and generalizes that of minimal ones. 

Moreover, we shall refer to biharmonic   submanifolds that are not minimal as \emph{proper} biharmonic submanifolds. 

Unfortunately, this notion of proper biharmonic   submanifolds appears too strong when the ambient space is the flat Euclidean space. In fact, for an immersion $\varphi:M^m\hookrightarrow \mathbb{R}^n$, the biharmonic condition \eqref{bitension-H} reduces to
\begin{align}
\label{biharmonic-Rn}
\Delta {\mathbf H}=(\Delta {\rm H}^1,\ldots,\Delta {\rm H}^n)=0
\end{align}
where $\Delta$ is the Beltrami-Laplace operator on $M$, and the mean curvature vector field is seen as a map ${\mathbf H}:M\to\mathbb{R}^n$.

By using the structure equations, B.-Y. Chen in \cite{MR1133117} and G.~Y.~Jiang in \cite{MR924896} proved that, in the case $m=2$ and  $n=3$, equation \eqref{biharmonic-Rn} implies that $|\mathbf H|^2$ is constant. Then, the vanishing of the component of \eqref{biharmonic-Rn} in the direction of the normal to the surface gives
\[
|B|^2 |\mathbf H|=0
\]
which implies that ${\mathbf H}=0$. In summary, we have that \emph{biharmonic   surfaces in $\mathbb{R}^3$ are necessarily minimal}.

This result, together with the lack of known examples in higher dimensions, motivated B.-Y. Chen to propose, in \cite{MR1143504}, the following conjecture, now widely known as Chen's Conjecture.\\

{\bf Conjecture 1}: {\em Any biharmonic   submanifold in $\mathbb{R}^n$ is minimal.}\\

Chen's conjecture is a local conjecture and has been verified in several cases. For instance, a few years after Chen's work, Hasanis and Vlachos proved in \cite{MR1330627} that it holds for hypersurfaces in $\mathbb{R}^4$. Furthermore, in \cite{MR4233281,MR4669301}, the conjecture has been confirmed for hypersurfaces in $\mathbb{R}^5$ and $\mathbb{R}^6$. However, the general conjecture remains open and in higher codimension, even the case of surfaces $M^2\hookrightarrow\mathbb{R}^4$ is not yet completely solved. 

Things do not change when considering biharmonic   submanifolds $\varphi:M^m\hookrightarrow (N^n,h)$ in a space of negative constant sectional curvature. In this case, the major challenge is to prove the following Generalized Chen's Conjecture.\\

{\bf Conjecture 2}: {\em Any biharmonic   submanifold in a space with non-positive constant sectional curvature is minimal.}\\

Also in this case there are several results supporting its validation (see, for example, \cite{MR2448058,MR3420554,MR3770116,MR2811621,MR2004799,MR2734173}). 

\begin{remark}\label{rmk-conjecture2}
We should point out that, according to a result of \cite{MR2004799}, to prove Conjecture 2 is enough to show that a biharmonic   submanifold in a space with non-positive sectional curvature has constant norm of the mean curvature vector field. 
\end{remark}

\begin{remark}
The original formulation of the Generalized Chen's Conjecture (see \cite{MR1871019}) assumes that the sectional curvature of the ambient space is non-positive but does not require it to be constant. In this generality, the conjecture is known to be false, as demonstrated by the existence of a non-minimal biharmonic hyperplane in $\R^5$ equipped with a conformally flat metric with non-constant negative curvature (see \cite{MR2975260}).
\end{remark} 

Fortunately, when the ambient space is positively curved, the situation change in positive. For instance, if we consider isometric immersions 
$\varphi:M^m\hookrightarrow \s^n$ to the $n$-dimensional round sphere, the biharmonic condition \eqref{bitension-H} becomes
\[
\bar\Delta\mathbf H=m\mathbf H
\]
and the following are main examples of proper solutions (see \cite{MR1863283,MR1919374,MR886529}):

\begin{itemize}
\item[{\bf B1}] The canonical inclusion of the small hypersphere
$$
\mathbb{S}^{n-1}(1 / \sqrt{2})=\left\{(x, 1 / \sqrt{2}) \in \mathbb{R}^{n+1}: x \in \mathbb{R}^n,|x|^2=1 / 2\right\} \subset \mathbb{S}^n.
$$    
\item[{\bf B2}] The canonical inclusion of the standard (extrinsic) products of spheres
$$
  \mathbb{S}^{n_1}(1 / \sqrt{2}) \times \mathbb{S}^{n_2}(1 / \sqrt{2})=\left\{(x, y) \in \mathbb{R}^{n_1+1} \times \mathbb{R}^{n_2+1},|x|^2=|y|^2=1 / 2\right\} \subset \mathbb{S}^n 
$$
$n_1+n_2=n-1 \text { and } n_1 \neq n_2$.
\item[{\bf B3}] The maps $\varphi=\imath \circ \phi: M \rightarrow \mathbb{S}^n$, where $\phi: M \hookrightarrow \mathbb{S}^{n-1}(1 / \sqrt{2})$ is a minimal immersion, and $\imath: \mathbb{S}^{n-1}(1 / \sqrt{2}) \hookrightarrow \mathbb{S}^n$ denotes the canonical inclusion.
\item[{\bf B4}] The maps $\varphi=\imath \circ\left(\phi_1 \times \phi_2\right): M_1 \times M_2 \hookrightarrow \mathbb{S}^n$, where $\phi_i: M_i^{m_i} \hookrightarrow \mathbb{S}^{n_i}(1 / \sqrt{2})$, $0<m_i \leq n_i, i=1,2$, are minimal immersions, $m_1 \neq m_2, n_1+n_2=n-1$, and $\imath: \mathbb{S}^{n_1}(1 / \sqrt{2}) \times \mathbb{S}^{n_2}(1 / \sqrt{2}) \hookrightarrow \mathbb{S}^n$ denotes the canonical inclusion.
\end{itemize}

If we restrict to the case of hypersurfaces $\varphi:M^m \hookrightarrow\s^{m+1}$, writting ${\mathbf H}=f \eta$ where $\eta$ is a unit normal vector field to $M$, the biharmonic condition \eqref{bitension-H} can be splited into its normal and tangential components to obtain that a hypersurfaces in the sphere is biharmonic if and only if
\begin{equation}\label{eq-biharhypersn}
\left\{\begin{array}{ll}
 \Delta f-\left(m-|A|^2\right) f =0&\text{Normal } \\
2A(\operatorname{grad} f)+m f \operatorname{grad} f=0 & \text {Tangential} 
\end{array}\right.
\end{equation}
where $A$ denotes the shape operator. According to the BMO conjecture (see \cite{MR2448058}), the main examples, {\bf B1} and {\bf B2}  should be the only two examples of proper biharmonic  hypersurfaces in $\s^{m+1}$. In other words we have the following\\

{\bf Conjecture 3}: {\em The examples {\bf B1} and {\bf B2} (with $n=m+1$) are the only  proper biharmonic hypersurfaces of the sphere $\s^{m+1}$.}\\

Positive solutions of Conjecture 3 are known in the following cases: for   surfaces in $\s^3$ \cite{MR1863283}; for compact hypersurfaces in $\s^4$ \cite{MR2560665}; for isoparametric hypersurfaces in $\s^{m+1}$, $m\geq 2$, \cite{MR2943022}. 

We observe that the main examples, {\bf B1} and {\bf B2}, represent CMC ($f={\rm constant}$) biharmonic hypersurfaces in the sphere $\s^{m+1}$. Also, from  \eqref{eq-biharhypersn}, we deduce immediately that a non-minimal CMC hypersurface $\varphi:M^m \hookrightarrow\s^{m+1}$ is proper biharmonic if and only if
\[
|A|^2=m.
\]

Thus Conjecture 3, for CMC proper biharmonic  hypersurfaces in $\s^{m+1}$, can be entirely rephrased as a conjecture in classical differential geometry as follows.\\

{\bf Conjecture 4}: {\em The examples {\bf B1} and {\bf B2} (with $n=m+1$) are the only non-minimal CMC hypersurfaces of the sphere $\s^{m+1}$ with $|A|^2=m$.}

\begin{remark} 
The analog of Conjecture 4 in the minimal case has a defenitive answer. This is a classical result proved by Lawson \cite{MR238229} and by Chern, do Carmo \& Kobayashi \cite{MR273546}: {\em minimal hypersurfaces in $\s^{m+1}$ with $|A|^2=m$ are products
$\s^p(a)\times\s^q(b)\hookrightarrow\s^{m+1}$ with $p+q=m$, $a^2=p/(p+q)$ and $b^2=q/(p+q)$.}
\end{remark}

\begin{remark}
Isoparametric hypersurfaces of the sphere $\s^{m+1}$ have constant principal curvatures, making them CMC with  $|A|^2={\rm constant}$. A classical problem in differential geometry, known as the Generalized Chern Conjecture,  aims to establish the converse: {\em a CMC hypersurface in $\s^{m+1}$ with $|A|^2={\rm constant}$ is isoparametric.} Since the only isoparametric hypersurfaces  in $\s^{m+1}$ with $|A|^2=m$ are the examples {\bf B1} and {\bf B2}, Conjecture 4 can  be regarded as a special case of the Generalized Chern Conjecture.
\end{remark}

Conjecture 3 and Conjecture 4 remain open and, consequently, the classification of  biharmonic hypersurfaces in the sphere $\mathbb{S}^{m+1}$ is still unresolved. \\

As a final remark, considering Conjecture 4,
the situation differs from that of hypersurfaces in spaces with non-positive sectional curvature (see Remark~\ref{rmk-conjecture2}). Specifically, for biharmonic hypersurfaces in $\s^{m+1}$ the condition $|{\mathbf H}|={\rm constant}$ does not automatically imply that is one of the known examples. 
However, to prove that a biharmonic hypersurface is CMC is not an easy task either.
Therefore, the following conjecture can be seen as a preliminary step -- though not a definitive one -- toward the classification of biharmonic hypersurfaces in $\s^{m+1}$.\\

{\bf Conjecture 5}: Any biharmonic hypersurface in $\s^{m+1}$ is CMC.\\

The latter Conjecture is actually a special case of the following more general conjecture.
\\

{\bf Conjecture 6}: Any biharmonic submanifold in $\s^{m+1}$ has $|{\mathbf H}|={\rm constant}$.\\

 For a deeper understanding of the state of the art regarding classification and rigidity results for biharmonic submanifolds in the sphere $\mathbb{S}^{m+1}$, a good starting point is the paper \cite{MR4410183} and the references therein. Additionally, \cite{MR4541401} provides valuable insights on this topic.

\section{Biconservative   submanifolds}

In this section we provide a different way to generalize the notion of minimal   submanifolds. The idea trace back to the work of Hilbert \cite{MR1512197} when he defined the notion of {\em stress-energy tensor} for a given variational problem. More precisely, the { stress-energy
tensor} associated to a variational problem is a symmetric $2$-covariant tensor $S$ conservative, that is  $\Div S=0$, at critical points of the corresponding functional.

In the context of harmonic maps as critical points of \eqref{energia}, the stress-energy tensor was studied in details by
Baird and Eells in~\cite{MR655417}. Indeed, they proved that the tensor  
\begin{equation}\label{strees-energy}
S(X,Y)=\frac{1}{2} |d\varphi|^2 g(X,Y)-\varphi^*h(X,Y)
\end{equation}
satisfies 
\[
\Div S(X)=-\langle\tau(\varphi),d\varphi(X)\rangle,
\]
 thus adhering to the principle of a stress-energy tensor, that is $\Div S=0$ when $\varphi$ is harmonic, i.e. when $\tau(\varphi)=0$.

If $\varphi:M^m\hookrightarrow (N^n,h)$ is an isometric immersion, we have already observed that $\tau(\varphi) = m{\mathbf H}$. Consequently, for a minimal   submanifold $\varphi:M^m\hookrightarrow (N^n,h)$, the stress-energy tensor is conservative. One might attempt to generalize the notion of minimal   submanifolds by considering those with conservative stress-energy tensor \eqref{strees-energy}. However, this approach yields no new family of immersions. Indeed, for any isometric immersion $\varphi:M^m\hookrightarrow (N^n,h)$, the tension field $\tau(\varphi) = m{\mathbf H}$ is normal to the immersed submanifold. As a result, we always have 
\[
\Div S = -m\langle{\mathbf H},d\varphi\rangle = 0,
\] 
regardless of the immersion. \\

Hoping for better luck, let us consider the stress-energy tensor associated to the bienergy \eqref{bienergy}. In this context, Jiang \cite{MR891928} (see also \cite{MR2395125}) constructed an ad-hoc $(0,2)$-tensor:  
\begin{equation}\label{eq-stress-bienergy}
\begin{aligned}
S_2(X,Y)=&\frac{1}{2}\vert\tau(\phi)\vert^2\langle X,Y\rangle+
\langle d\phi,\nabla\tau(\phi)\rangle \langle X,Y\rangle \\
 & -\langle d\phi(X), \nabla_Y\tau(\phi)\rangle-\langle
d\phi(Y), \nabla_X\tau(\phi)\rangle,
\end{aligned}
\end{equation}
which satisfies the relation $\Div S_2=-\langle\tau_2(\phi),d\phi\rangle$, thus
conforming to the principle of a  stress-energy tensor for the bienergy: $\Div S_2=0$ when $\varphi$ is a biharmonic map, i.e. when  $\tau_2(\phi)=0$. We shall call $S_2$ the {\em stress-bienergy} tensor.

\begin{remark}
It is interesting to note that both the stress-energy tensor and the stress-bienergy tensor can be characterized through metric variations rather than variations of the map. For the harmonic case, see the paper by Sanini \cite{MR710809}, and for the biharmonic case, see \cite{MR2395125}. 
\end{remark}

Now, we restrict our focus to isometric immersions $\varphi:M^m\hookrightarrow (N^n,h)$ and revisit the question we previously posed for the stress-energy tensor $S$:  \\

{\bf Question 1}: {\em Can we study the isometric immersions satisfying $\Div S_2=0$?} \vspace{3mm}

In this case, the answer is affirmative, since the bitension field $\tau_2(\varphi)$ of a submanifold is not always normal to the immersion. Consequently, the condition  
\[
\Div S_2=-\langle \tau_2(\varphi), d\varphi\rangle =0
\]
clearly identifies a new class of   submanifolds. This reasoning led to the following definition, introduced in \cite{MR3180932}.

 \begin{definition}
 An immersed submaifold $\varphi:M\hookrightarrow(N,h)$ is called {\em biconservative} if $\Div S_2=0$.
 \end{definition}
 
 We have the following direct consequences.

\begin{enumerate}

\item Any minimal   submanifold is also biconservative.

\item A   submanifold $\varphi:M\hookrightarrow(N,h)$  is biconservative if and only if the tangential component of the bitension field is identically zero, that is $\tau_2(\varphi)^{\top}=0$.
\end{enumerate}

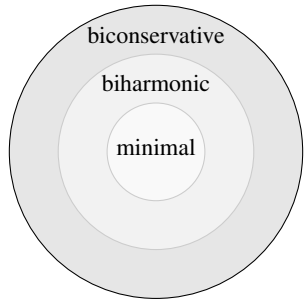
\begin{figure}[h!]
\centering
\begin{tikzpicture}[scale=0.65]
    \fill[gray!20] (0,0) circle[radius=3];
    \draw[black] (0,0) circle[radius=3];
    \node at (0,2.4) {\scriptsize \textcolor{black}{biconservative}};

    \fill[gray!10] (0,0) circle[radius=2];
    \draw[black, gray!40] (0,0) circle[radius=2];
    \node at (0,1.4) {\scriptsize \textcolor{black}{biharmonic}};

    \fill[gray!5] (0,0) circle[radius=1];
    \draw[black, gray!40] (0,0) circle[radius=1];
    \node at (0,0.1) {\scriptsize minimal};
\end{tikzpicture}
\caption{The inclusions between minimal, biharmonic, and biconservative   submanifolds.}
\label{fig:inclusions1}
\end{figure}

We thus have the inclusions between the families of minimal, biharmonic and biconservative   submanifolds, see Figure~\ref{fig:inclusions1}.

Nowadays,  the theory of biconservative   submanifolds has become a growing area of study, attracting numerous contributions to the field. For an overview, the reader may refer to the article by Fetcu \& Oniciuc \cite{MR4410183} or to the survey by B.-Y. Chen \cite{chen2024recentdevelopmentbiconservativesubmanifolds}.

In this note, instead of presenting a survey of the current state of the art in the field,  we would like just  to highlight some properties of biconservative immersions that could justify their studies. The focus will be on surfaces $\varphi:M^2\hookrightarrow (N^n,h)$.

\subsection{Biconservative surfaces and the generaliszed Hopf function}
For an immersion $\varphi:M^2\hookrightarrow N^3(c)$ of an oriented surface in a three-dimensional space form of constant sectional curvature $c$,
denote by $g=\langle , \rangle$ the induced metric on $M^2$. By assumption, $M^2$ is orientable and then it is a one-dimensional complex manifold. If we consider local isothermal coordinates
$(U; x,y)$, then  $g= \lambda^2 ( dx^2 + dy^2)$ for some positive function $\lambda$ on $U$ and $\{\partial_x,\partial_y\}$
is positively oriented. Let us denote, as usual,
$$
z=x+\i y\,,\quad
\partial_z=\frac{\partial}{\partial z}=\frac{1}{2}\left(\frac{\partial}{\partial x}-\i \frac{\partial}{\partial y}\right)\,,\quad
\partial_{\bar{z}}=\frac{\partial}{\partial \bar{z}}=\frac{1}{2}\left(\frac{\partial}{\partial x}+\i \frac{\partial}{\partial y}\right)\,.
$$
Then, the classical Hopf function is defined by
\begin{equation}\label{eq:hopf-classical}
\Phi(z,\bar{z})=\langle B(\partial_z,\partial_z),\eta \rangle=\langle A(\partial_z),\partial_z \rangle.
\end{equation}

The Hopf function defined in \eqref{eq:hopf-classical} is the key ingredient in the proof of the famous Hopf's Theorem: {\it a CMC immersed sphere in $N^3(c)$ is a round sphere}. 

Among others, Hopf's proof is based on the following fact: {\em $\Phi$ is holomorphic if and only if $\varphi:M^2\hookrightarrow N^3(c)$ is CMC}. 

The proof of the latter follows by a straightforward computation, taking into account Codazzi's equation, which gives
\begin{equation}\label{eq-phi-olomorphic}
2\,\partial_{\bar{z}}\Phi(z,\bar{z})=\lambda^2 \partial_{{z}}(f)\,.
\end{equation} 

In higher dimension, that is for immersions $\varphi:M^2\hookrightarrow N^n(c)$, $n\geq 4$, the function $\Phi$ cannot be defined and a fair substitute is the function
\begin{equation}\label{eq:Q-classical}
Q(z,\bar{z})=\langle B(\partial_z,\partial_z),{\mathbf H}\rangle=\langle A_{\mathbf H}(\partial_z),\partial_z \rangle\,,
\end{equation}
where $A_{\mathbf H}$ is the shape operator in the direction of the normal vector field ${\mathbf H}$.
 A classical result, see \cite{MR390973, MR370443}, states that {\em if the surface $\varphi:M^2\hookrightarrow N^n(c)$, $n\geq 4$, has parallel mean curvature (PMC) then $Q$ is holomorphic}.
 
A natural question is to ask whether the converse holds, that is:\\

{\bf Question 2}: {\em When a surface $\varphi:M^2\hookrightarrow N^n(c)$, $n\geq 4$, with $Q$ holomorphic  is PMC or CMC?}\\

 This question
is interesting also for surfaces $M^2\hookrightarrow N^3(c)$. In fact,  taking into account \eqref{eq-phi-olomorphic}, we easily obtain
\begin{equation}\label{eq-Q-dzbar}
4\,\partial_{\bar{z}}Q(z,\bar{z})= \lambda^2 \left[  \langle A(\grad f),\partial_x \rangle-  \i  \langle A(\grad f),\partial_y \rangle\right].
\end{equation}

As a consequence of \eqref{eq-Q-dzbar}, if we denote by $K^M$ the Gaussian curvature of the surface $M$, we obtain
\begin{proposition}\label{pro:Q-holomorphic}
Let $M^2\hookrightarrow N^3(c)$ be an oriented surface in a space form of constant sectional curvature $c$.
\begin{itemize}
\item[(a)] If $\det(A)=K^M-c\neq 0$, then $f$ is constant if and only if $Q(z,\bar{z})$ is holomorphic;
\item[(b)] If $f$ is not constant and  $Q(z,\bar{z})$ is holomorphic, then $K^M=c$.
\end{itemize}
\end{proposition}
We point out that surfaces satisfying condition (b) of Proposition~\ref{pro:Q-holomorphic} do exist. For instance, in $\R^3$ 
the cone $(z-1)^2 =x^2 + y^2$, $z>0,$ has $Q$ holomorphic but it is not CMC. All surfaces satisfying condition (b) of Proposition~\ref{pro:Q-holomorphic} can be actually classified, see \cite{MR3441516} for details.\\

Going back to Question 2 we shall now give an answer to it without the assumption that the ambient space is of constant sectional curvature. More precisely, we shall consider: 
 \\

{\bf Question 3}: {\em Does there exist a class of immersed surfaces $\varphi:M^2\hookrightarrow N^n$ for which the holomorphicity of 
$Q$ is equivalent to being CMC?}\\

Surprisingly the answer of Question 3 is intimately related to biconservative surfaces as we shall show below. We begin with a rather general fact

\begin{proposition}[see \cite{MR3441516} and \cite{MR3693945} for a generalized version]\label{pro:T-ho-codazzi}
Let $T$ be a symmetric $(0,2)$-tensor field on a Riemannian surface $(M^2,g)$ and set $t=\trace T$. Assume that $M^2$ is orientable and $\Div T=0$ on $M^2$. Then {$T(\partial_z,\partial_z)$ is holomorphic if and only if $t=\cst$}.
\end{proposition}

Let now $M^2\hookrightarrow N^n$ be a surface in an $n$-dimensional Riemannian manifold and denote by $g$ the induced metric. Then the stress energy tensor $S_2$, defined in \eqref{eq-stress-bienergy}, is indeed a symmetric $(0,2)$-tensor
on $(M^2,g)$ and, in this case,  the expression \eqref{eq-stress-bienergy} of $S_2$ reduces to
$$
{S_2(X,Y)=-2  |{\mathbf H}|^2 \langle X, Y\rangle + 4 \langle A_{\mathbf H}(X), Y\rangle}\,.
$$
Taking the trace of $S_2$ we have
$$
t=\trace S_2=4 |{\mathbf H}|^2.
$$
Moreover, since $g(\partial_z,\partial_z)=0$,
we obtain
$$
S_2(\partial_z,\partial_z)=4\,\langle A_{\mathbf H}(\partial_z),\partial_z\rangle=4\, Q(z,\bar{z})\,.
$$

\begin{figure}[h!]
\centering
\begin{tikzpicture}[scale=0.30]
\fill[gray!20] (0,0) ellipse (8 and 5);
    \node at (0,3.8) {\scriptsize \textcolor{black}{In the set of}};
    \node at (0,3) {\scriptsize \textcolor{black}{biconservative surfaces}};
    \node at (0,2.1) {\scriptsize $M^2\hookrightarrow N^n$};

    \node at (0,0) {\bf \scriptsize \textcolor{black}{CMC \; $\Leftrightarrow$ \; $Q$ is holomorphic}};
\end{tikzpicture}
\caption{The solution to Question 3.}
\label{fig:biconservativeQholo}
\end{figure}
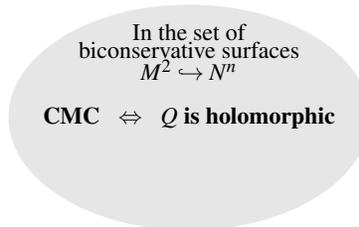

Thus, as a direct consequence of Proposition~\ref{pro:T-ho-codazzi}, we obtain the following answer to Question 3 which offers a compelling justification for studying biconservative surfaces (see also Figure~\ref{fig:biconservativeQholo}).

\begin{theorem}\label{teo:Q-hol-h-constant-biconservative}
Let $M^2\hookrightarrow (N^n,h)$ be a biconservative surface in an $n$-dimensional Riemannian manifold. Then
$Q(z,\bar{z})$ is holomorphic if and only if $|{\mathbf H}|=\cst$.
\end{theorem}

\subsection{Biconservative surfaces in $3$-dimensional space forms} In this last part we shall concentrate to the case of   surfaces $M^2$ in a $3$-dimensional space form $N^3(c)$ of constant sectional curvature $c$. As already pointed out, a  surface $\varphi:M^2\hookrightarrow N^3(c)$ is biconservative if and only if the tangential component of the bitension field is identically zero, that is,  taking into account that the tangential component of $\tau_2$ given in \eqref{eq-biharhypersn} remains the same when the ambient space is replaced by any space with constant sectional curvature,  if and only if
\begin{equation}\label{eq:bi-suf-3sp}
A(\grad f) + f \grad f = 0.
\end{equation}

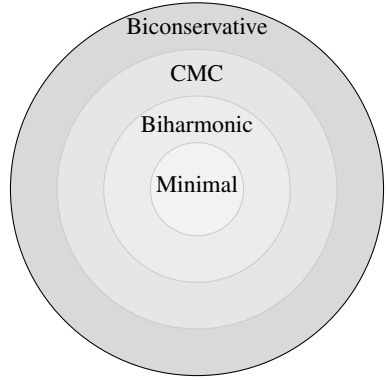
\begin{figure}[h]
\centering
\begin{tikzpicture}[scale=0.62]
    \fill[gray!30] (0,0) circle[radius=4cm];
    \draw[black] (0,0) circle[radius=4cm];
    \node at (0,3.5) {\scriptsize Biconservative};

    \fill[gray!20] (0,0) circle[radius=3];
    \draw[black,gray!40] (0,0) circle[radius=3];
    \node at (0,2.5) {\scriptsize CMC};

    \fill[gray!15] (0,0) circle[radius=2];
    \draw[black,gray!40] (0,0) circle[radius=2];
    \node at (0,1.4) {\scriptsize Biharmonic};

    \fill[gray!10] (0,0) circle[radius=1];
    \draw[black,gray!40] (0,0) circle[radius=1];
    \node at (0,0.1) {\scriptsize Minimal};
\end{tikzpicture}
\caption{Inclusions between the families of minimal, biharmonic, CMC, and biconservative surfaces in a $3$-dimensional space form $N^3(c)$.}
\label{fig:inclusion3}
\end{figure}

We then obtain immediately from \eqref{eq:bi-suf-3sp} that a CMC surface in a $3$-dimensional space form is biconservative. Moreover, we have already mentioned that B.-Y. Chen proved in \cite{MR1143504} that a biharmonic surface in $\R^3$ is CMC. The same result was also proved in \cite{MR1863283,MR1919374} for surfaces in $N^3(c)$ with $c\neq 0$. We then have the inclusions shown in Figure~\ref{fig:inclusion3} between minimal, biharmonic, CMC and biconservative surfaces in $N^3(c)$.

From this, it becomes clear that the primary interest is on non-CMC  biconservative surfaces, that is $\grad f\neq 0$ at some points. Let   $\varphi:M^2\hookrightarrow N^3(c)$ be a biconservative surface and assume, for simplicity, that $\grad f\neq 0$ everywhere on $M$. Then, from \eqref{eq:bi-suf-3sp}, the vector field
\[
e_1=\frac{\grad f}{|\grad f|}
\]
is a principal direction with corresponding principal curvature 
\[
\lambda_1=-f=-\frac{1}{2}(\lambda_1+\lambda_2)
\] 
where $\lambda_2$ is the principal curvature corresponding to the principal direction $e_2\perp e_1$. Consequently, biconservative surfaces $\varphi:M^2\hookrightarrow N^3(c)$ are {\em Linear Weingarten} surfaces with
\[
3\lambda_1+\lambda_2=0.
\]
This served as the starting point for obtaining a geometric characterization of biconservative surfaces $\varphi:M^2\hookrightarrow N^3(c)$, as stated in the following result, originally proved in \cite{MR3180932} and later presented in this form in \cite{MR4685022}.  
\begin{proposition}
Let $M^{2}$ be a surface of a space form $N^3(c)$ with nowhere zero $\grad f$. Then, $M^{2}$ is {biconservative} if and only if it is {rotational} and the principal curvatures satisfy {$3\lambda_1+\lambda_2=0$}.
Moreover, the profile curve lies in a totally geodesic surface $N^2(c)\subset N^3(c)$ and its curvature is given by {$\kappa=-\lambda_1$}.
\end{proposition}

\begin{example}
For example, a biconservative  surface $\varphi:M^2\hookrightarrow \R^3$ with $\grad f\neq 0$ at {any point} and $f > 0$  is, locally,  the $SO(2)$-invariant immersion  

$$
\varphi(\rho,\vartheta)=(x(\rho) \cos \vartheta, x(\rho) \sin \vartheta, z(\rho))
$$
with
$$
\frac{dz}{dx}=\frac{C\, x^{1/3}}{\sqrt{x^{2/3}-C^2}}\,,\quad C>0\,.
$$

In the Figure~\ref{fig:profilebiconservativer3}, the continuous line represents a plot of the profile curve of the local surface described above. The figure also includes the dashed symmetric extension of the profile curve. Together, the continuous and dashed lines define the complete profile curve of a complete $SO(2)$-invariant biconservative surface in $\mathbb{R}^3$.  

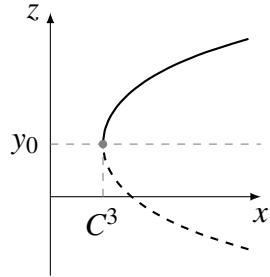
\begin{figure}[h]
\centering
\begin{tikzpicture}[xscale=0.7, yscale=0.7]
    \draw [->] (0,0) -- (4,0) node[below] {$x$};
    \draw [->] (0,-1.5) -- (0,3.5) node[left] {$z$};

    \draw [thick] plot [smooth] coordinates { 
        (1., 1.) (1.03141, 1.25) (1.12763, 1.5) (1.29468, 1.75) 
        (1.54308, 2.) (1.88842, 2.25) (2.35241, 2.5) (2.96419, 2.75) 
        (3.7622, 3.) 
    };

    \draw [gray, dashed] (1,0) -- (1,1);
    \draw [gray, dashed] (0,1) -- (4,1);

    \draw (1,0) node[below] {$C^3$};
    \draw (0,1) node[left] {$y_0$};

    \draw [dashed, thick] plot [smooth] coordinates { 
        (3.7622, -1.) (2.96419, -0.75) (2.35241, -0.5) (1.88842, -0.25) 
        (1.54308, 0.) (1.29468, 0.25) (1.12763, 0.5) (1.03141, 0.75) (1., 1.) 
    };

    \filldraw [gray] (1,1) circle (2pt);
\end{tikzpicture}
\caption{Plot of the profile curve of a biconservtive surface in $\R^3$.}
\label{fig:profilebiconservativer3}
\end{figure}
\end{example}
Gluing, along their boundaries, maximal biconservative surfaces with $\grad f\neq 0$ at any point, we obtain smooth non-CMC biconservative surfaces in $N^3(c)$, with $\grad f\neq 0$ at any point of an open dense subset of the resulting domain. The behavior of such surfaces, particularly concerning their completeness, was studied in detail by Nistor in \cite{MR3566105} and by Nistor \& Oniciuc in \cite{MR4553530}. In \cite{MR3566105}, Nistor posed the question of whether, among the surfaces constructed above,  {\em closed} (compact without boundary) biconservative surfaces exist in $\mathbb{S}^3$. In the final part of this excursion, we present an approach to finding a positive answer to this question.
First recall the following interpretation of biconservative surfaces in space forms.
\begin{proposition}[\cite{MR4685022}]
A  biconservative surface in a $3$-dimensional space form \( N^3(c) \), with nowhere zero $\grad f$,  is a rotational surface whose profile curve is a critical point of the following {\em bending-energy} functional:
\[
{\Theta}(\gamma) := \int_\gamma \kappa^{1/4} \, ds \quad \text{defined on the space of curves in } N^2(c).
\]
Thus, the curvature \(\kappa\) of the profile curve satisfies the corresponding Euler-Lagrange equation:
\begin{equation}\label{eq:fund-kappa}
\kappa^{3/4} \frac{d^2}{ds^2} \left( \frac{1}{\kappa^{3/4}} \right) - 3\kappa^2 + c = 0.
\end{equation}
The converse is also true.
\end{proposition}

Now \eqref{eq:fund-kappa} admits the following prime integral
\begin{equation}\label{eq:fund-kappa2}
\kappa_s^2=\frac{16}{9}\kappa^2\left(16\, d \,\kappa^{3/2}-9\,\kappa^2-c\right),\qquad d\in\R,
\end{equation}
which can be written, putting $u=\sqrt{\kappa}$, as
\begin{equation}\label{eq:fund-kappa3}
u_s^2=\frac{4}{9}u^2\left(16\, d\, u^3-9\, u^4-c\right)=\frac{4}{9} u^2 Q(u).
\end{equation}

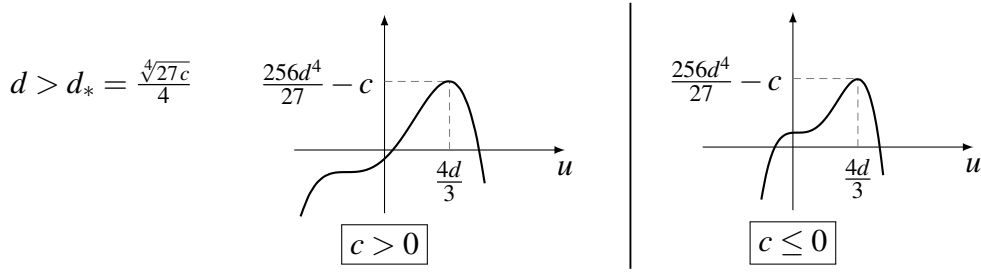
\begin{figure}[h!]
\centering
\begin{tikzpicture}[scale=1.2]
    \draw[->] (-1,0) -- (2,0) node[below] {${\small u}$};
    \draw[->] (0,-.7) -- (0,1.5);
    
    \draw[gray,densely dashed] (.72,0)--(.72,.76);
    \draw[gray,densely dashed] (0,.76)--(.72,.76);
    
    \draw (.72,0) node[below] {${\small \frac{4d}{3}}$};
    \draw (0,.76) node[left] {${\small \frac{256 d^4}{27}-c}$};
    \draw (-2,.76) node[left] {${\small d>d_*=\frac{\sqrt[4]{27\,c}}{4}}$};
    \draw (0,-.7) node[below] {$\boxed{c> 0}$};

    \draw[thick, xshift=-.42cm, yshift=-.075cm, scale=1.7] 
    plot [smooth] coordinates {   
        (-0.3, -0.3889) (-0.25, -0.260156) (-0.2, -0.1784) 
        (-0.15, -0.131556) (-0.1, -0.1089) (-0.05, -0.101056) (0., -0.1) 
        (0.05, -0.0990563) (0.1, -0.0929) (0.15, -0.0775563) (0.2, -0.0504) 
        (0.25, -0.0101563) (0.3, 0.0431) (0.35, 0.107944) (0.4, 0.1816) 
        (0.45, 0.259944) (0.5, 0.3375) (0.55, 0.407444) (0.6, 0.4616) 
        (0.65, 0.490444) (0.7, 0.4831) (0.75, 0.427344) (0.8, 0.3096) 
        (0.85, 0.114944) (0.9, -0.1729) 
    };
\end{tikzpicture}
\hspace{4mm}
\begin{tikzpicture}[scale=1.2]
    \draw[black,line width=0.25mm] (-1.8,-1.35) -- (-1.8,1.6);

    \draw[->] (-1,0) -- (2,0) node[below] {${\small u}$};
    \draw[->] (0,-.7) -- (0,1.5);
    
    \draw[gray,densely dashed] (.72,0)--(.72,.76);
    \draw[gray,densely dashed] (0,.76)--(.72,.76);
    
    \draw (.72,0) node[below] {${\small \frac{4d}{3}}$};
    \draw (0,.76) node[left] {${\small \frac{256 d^4}{27}-c}$};
    \draw (0,-.7) node[below] {$\boxed{c\leq 0}$};

    \draw[thick, xshift=.05cm, yshift=.26cm] 
    plot [smooth] coordinates { 
        (-0.4, -0.8424) (-0.35, -0.578056) (-0.3, -0.3889) (-0.25, -0.260156) 
        (-0.2, -0.1784) (-0.15, -0.131556) (-0.1, -0.1089) (-0.05, -0.101056) 
        (0., -0.1) (0.05, -0.0990563) (0.1, -0.0929) (0.15, -0.0775563) 
        (0.2, -0.0504) (0.25, -0.0101563) (0.3, 0.0431) (0.35, 0.107944) 
        (0.4, 0.1816) (0.45, 0.259944) (0.5, 0.3375) (0.55, 0.407444) 
        (0.6, 0.4616) (0.65, 0.490444) (0.7, 0.4831) (0.75, 0.427344) 
        (0.8, 0.3096) (0.85, 0.114944) (0.9, -0.1729) (0.95, -0.571556) 
    };
\end{tikzpicture}
\caption{The graph of $Q(u)$ according to the value of the curvature $c$.}
\label{fig:plotofQ}
\end{figure}

The graph of $Q(u)$, depending on the value of the curvature $c$, is shown in Figure~\ref{fig:plotofQ}. Note that only the region with $u > 0$ should be considered.

Setting $u=x^2$ and $y=x_s$ \eqref{eq:fund-kappa3} becomes\
\[
 y^2=\frac{x^2}{9} Q(x^2)
\]
which represents an algebraic curve $C$. A standard analysis, using the square root argument, shows that the trace of the curve $C$ is closed when $c>0$ and $ d>d_*=\sqrt[4]{27\,c}/4$, while it is not closed when $c\leq 0$.
Thus, if $c>0$ the curve $\alpha(s)=(x(s), y(s))$ is a curve in the trace of the closed curve (C).
When the curve $\alpha$ is defined in the maximal interval, since it can be seen as the integral curve of a smooth vector field without singularities, from the {\em Poincaré-Bendixon Theorem}, we deduce that $\alpha(s)$ is periodic. Consequently, $x(s)$, $u(s)=x^2(s)$  and $\kappa(s)=u(s)^2$ are periodic. \\

In summary, by setting $c=1$, we have shown that the profile curve of a non-CMC biconservative surface in $\mathbb{S}^3$ has periodic curvature. Since the surface is rotational, to guarantee the existence of a closed biconservative surface in $\mathbb{S}^3$, we must demonstrate that the profile curve $\gamma$ itself is periodic, not just its curvature.  

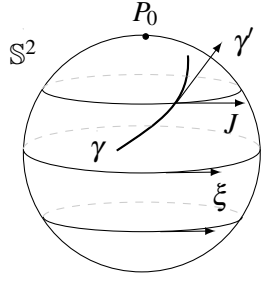
\begin{figure}[h]
\begin{center}
\begin{tikzpicture}[x=0.75pt,y=0.75pt,yscale=-.7,xscale=.7]

\draw   (364.67,150.67) .. controls (364.67,103.54) and (402.87,65.33) .. (450,65.33) .. controls (497.13,65.33) and (535.33,103.54) .. (535.33,150.67) .. controls (535.33,197.79) and (497.13,236) .. (450,236) .. controls (402.87,236) and (364.67,197.79) .. (364.67,150.67) -- cycle ;

\draw (364.67,60.67) circle(.01) node[below] {{\small \(\mathbf \s^2\)}};

\draw [color=black, thick]    (432,148.5) .. controls (463.25,127.75) and (485.25,114) .. (483.25,80) ;

\draw (417.75,134.25) circle(.01) node[below] {{\small \(\gamma\)}};

\draw[->]    (474.83,114.25) -- (524.75,114.25) ;

\draw[fill] (452.83,66.25) circle(2) node[above] {{\small $P_0$}};

\draw[color=black!20,dashed]   (364.67,147.67) .. controls (364.67,138.28) and (402.87,130.67) .. (450,130.67) .. controls (497.13,130.67) and (535.33,138.28) .. (535.33,147.67) ;

\draw[color=black] (535.33,147.67).. controls (535.33,157.06) and (497.13,164.67) .. (450,164.67) .. controls (402.87,164.67) and (364.67,157.06) .. (364.67,147.67) ;

\draw [color=black!20,dashed]   (378.33,104.33) .. controls (378.33,98.44) and (410.42,93.67) .. (450,93.67) .. controls (489.58,93.67) and (521.67,98.44) .. (521.67,104.33)  ;
\draw [color=black] (521.67,104.33).. controls (521.67,110.22) and (489.58,115) .. (450,115) .. controls (410.42,115) and (378.33,110.22) .. (378.33,104.33)  ;


\draw [color=black!20,dashed]   (379,196.33) .. controls (379,190.44) and (411.09,185.67) .. (450.67,185.67) .. controls (490.25,185.67) and (522.33,190.44) .. (522.33,196.33) ;

\draw [color=black]    (432,148.5) .. controls (463.25,127.75) and (485.25,114) .. (483.25,80) ;
\draw [color=black] (522.33,196.33).. controls (522.33,202.22) and (490.25,207) .. (450.67,207) .. controls (411.09,207) and (379,202.22) .. (379,196.33) ;
\draw[->]    (467.33,164.5) -- (507.25,164) ;
\draw (507.25,164) circle(.01) node[below] {{\small \(\xi\)}};
\draw[->]    (463.33,207) -- (503.25,206.5) ;

\draw [->]   (474.83,114.25) -- (507.83,70.25) ;
\draw (507.83,70.25) circle(.01) node[right] {{\small \(\gamma'\)}};

%
%
%
%
%
\draw (514.75,114.25) circle(.01) node[below] {{\small \(J\)}};
\end{tikzpicture}
\end{center}
\caption{The profile curve $\gamma$ in $\s^2$,  the integral curves of the Killing vector field $\xi$ and its restriction $J$ along $\gamma$.}
\label{fig:gammainsfere}
\end{figure}

To this end, let $\gamma: I \to \mathbb{S}^2$
be the profile curve of a biconservative surface \( M^2 \hookrightarrow \mathbb{S}^3 \). We can verify, using \eqref{eq:fund-kappa}, that the vector field
\[
J = -\frac{1}{4} \kappa^{1/4} T - \frac{3}{4} \kappa^{-7/4} N\vspace{-1mm}
\]
is a Killing vector field along \(\gamma\) in the sense of Langer \& Singer \cite{MR772124}. Then, according to \cite{MR772124}, \( J \) is the restriction to \(\gamma\) of a Killing vector field \(\xi\) of \(\mathbb{S}^2\).  Therefore, we can choose spherical coordinates $x(\vartheta,\psi)= (\cos\vartheta \sin\psi, \sin\vartheta \sin\psi, \cos\psi)$ of $\s^2$ so that its equator gives the only integral
geodesic of $\xi$. We shall denote by $P_0$ the north pole in the chosen spherical coordinates (see Figure~\ref{fig:gammainsfere}).

At this point, denoting by $\rho$ the period of the curvature $\kappa$, a standard argument (see, for example, \cite{MR2007599}), establishes that the curve $\gamma$ is periodic if there exist two integers $m,n\in\z$, with no common factor, such that its progression angle in one period of the curvature satisfies
\[
\int_0^{\rho} \vartheta'(s)\, ds= \frac{n}{m} 2\pi\,,
\]
where:
\begin{itemize}
\item the integer $m$ indicates that the period of $\gamma$ is $T=m \rho$;
\item the integer $n$ indicates how many times the curve goes around the pole $P_0$ of  $\s^2$.
\end{itemize}

Finally, using  the prime integral \eqref{eq:fund-kappa2}, we compute 
\[
I(d)=\int_0^{\rho} \vartheta'(s)\, ds=\int_0^{\rho} \frac{12 \kappa^{7/4}\sqrt{d}}{16 d \kappa^{3/2}-1}\, ds
\]
and obtain the following technical result, proven in \cite{MR4685022}.

\begin{proposition}\label{pro-technical}
The function $I(d)$ is strictly decreasing in $d$ and satisfies, for any $d\in\left( d_*,+\infty\right)$, 
\[
\pi<I(d)<\sqrt{2}\pi .
\]
\end{proposition}

By selecting integers $m$ and $n$ such that 
\[
\gcd(m, n) = 1 \quad \text{and} \quad m < 2n < \sqrt{2} m
\]
then we have
\[
\pi < \frac{n}{m} \cdot 2\pi < \sqrt{2} \pi.
\]
Thus, for these chosen values of $m$ and $n$, Proposition~\ref{pro-technical}  guarantees the existence of a $d\in\left( d_*,+\infty\right)$ such that:
\[
I(d) = \frac{n}{m} \cdot 2\pi.
\] 
In conclusion, we establish the announced existence result.
\begin{theorem}\label{teo-min-periodic}
There exists a discrete, biparametric family of closed non-CMC biconservative surfaces in the round 3-sphere $\mathbb{S}^3$.
\end{theorem}

\begin{remark}
None of the surfaces in Theorem~\ref{teo-min-periodic} is {\em embedded} in $\mathbb{S}^3$. In fact, it is not hard to deduce that for the surface to be embedded, the profile curve must close in a single round. That is, when $n=1$, there must exist an integer $m$ such that:
\[
\pi < \frac{2\pi}{m} < \sqrt{2}\pi
\]
which is not possible.
\end{remark}
\begin{example}
Choosing, as an example,  $m=3$ and $n=2$ we can (numerically) plot both the profile curve in $\s^2$ of the non-CMC  closed biconservative surface and its stereographic projection in $\R^3$, see Figure~\ref{fig:closedbiconservative}.
\begin{figure}[h]
\centerline{\includegraphics[width=.75\textwidth]{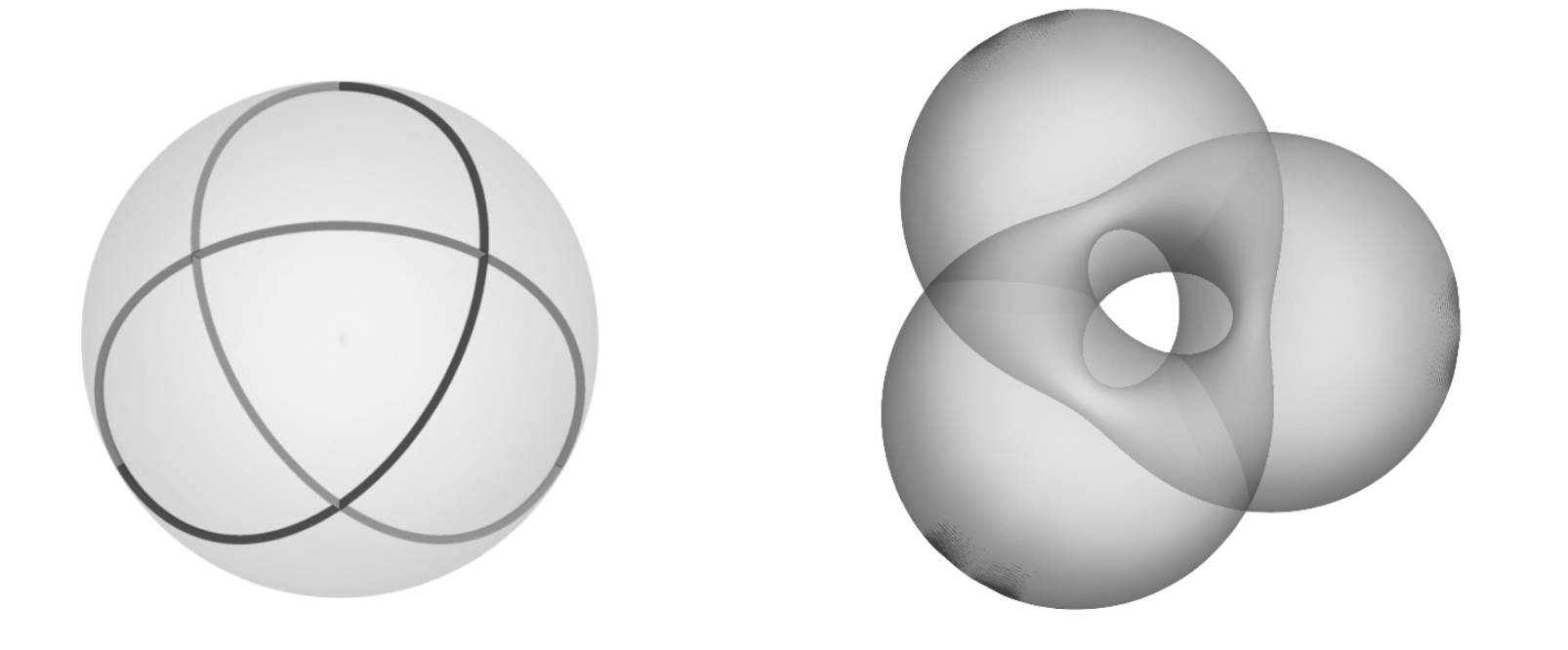}}
\caption{The profile curve (on the left) and the stereographic projection of the corresponding closed biconservative surface (on the right) in the case $m=3$ and $n=2$.}\label{fig:closedbiconservative}
\end{figure}
\end{example}

\begin{remark}
We point out that the full classification of complete, simply connected or not, non-CMC biconservative surfaces in $N^3(c)$ was finally obtained in \cite{MR4553530}.
\end{remark}

\begin{remark}
The construction of closed biconservative surfaces in $\s^3$ was extended in \cite{MR4489868}, demonstrating the existence of a discrete two-parameter family of non-CMC closed biconservative hypersurfaces in $\s^{m+1}$ for any 
$m\geq 2$.
\end{remark}

The theory of biconservative surfaces can be naturally described in three dimensional homogeneous spaces. In this context, it would be natural to pursuit the following problem.\\

{\bf Open problem}. Study the compact non-CMC biconservative surfaces in three dimensional homogeneous spaces.\\

{\bf Acknowledgments}. The author wishes to express sincere gratitude to the organizers of the Second Conference on Differential Geometry, held in Fez in October 2024, for their warm hospitality and meticulous organization of the event. The author also thanks Cezar Oniciuc and Andrea Ratto for their valuable feedback on an initial draft of this paper.


   \end{document}